\documentclass[10pt,reqno]{amsart}
\usepackage{geometry}
\usepackage{mathtools}
\usepackage{amsfonts,amsmath,amsthm,amssymb}
\usepackage{color}
\usepackage{hyperref,easybmat}

\numberwithin{equation}{section}

\newtheorem{thm}{Theorem}[section]
\newtheorem{openprob}[thm]{Open Problem}
\newtheorem*{openprob*}{Open Problem}
\newtheorem*{AC}{Alekseevskii Conjecture}

\newtheorem*{thm*}{Theorem}

\newtheorem{question}[thm]{Question}
\newtheorem*{question*}{Question}

\newtheorem*{observ*}{Observation}

\newtheorem{example}[thm]{Example}
\newtheorem{remark}[thm]{Remark}
\newtheorem*{remark*}{Remark}

\newcommand{\ip}[1]{\langle #1 \rangle}
\newcommand{\innerprod}{\ip {\cdot , \cdot} }

\begin{document}

\title{Survey: homogeneous Einstein manifolds}
\author{M.Jablonski}
\address{Department of Mathematics\\
	University of Oklahoma\\
	Norman, OK 73019-3103}
\email{mjablonski@math.ou.edu}
\date{November 17, 2021}
\maketitle

\setcounter{tocdepth}{1}
\tableofcontents

In Riemannian geometry, we have three standard notions of curvature which give a glimpse at the shape of a given metric and help us to discern how a space might be related to other familiar, model spaces.  These three curvatures are sectional, Ricci, and scalar curvature.

The spaces of constant sectional curvature are well-understood with their simply-connected covers being spheres, hyperbolic spaces, and Euclidean spaces.  
At the other extreme is scalar curvature.  All manifolds admit metrics of constant negative scalar curvature \cite{Lohkamp:MetricsOfNegativeRicciCurvature}.  In the homogeneous setting, scalar curvature is always constant and reveals some interesting, but limited, information about the underlying space \cite{Berard-Bergery:SurLaCourbureDesMetriquesRiemanniennesInvariantesDesGroupesDeLieEtDesEspacesHomegenes}.

Perhaps in the Goldilocks zone is  Ricci curvature.  If one were to hope to endow a given manifold with a special geometry against which other metrics could be weighed, this seems to be a reasonable invariant to hold constant and study \cite{Besse:EinsteinMflds}.  A Riemannian metric $g$ is called Einstein if it has constant Ricci curvature, i.e.
	\begin{align}
	\label{eqn: Einstein definition}
	ric_g = cg,
	\end{align}	 
for some $c\in\mathbb R$.

The general setting of manifolds  differs substansially from the homogeneous setting in both technique and results. For the reader interested in the general setting, we direct them to the survey \cite{Anderson:ASurveyOfEinsteinMetricsOn4Manifolds}.  Our interest is solely in the homogeneous setting and follows on two excellent surveys,  written almost a decade ago;  see  \cite{Wang:EinsteinMetricsFromSymmetryAndBundleConstractionsASurvey}  for the compact setting and \cite{Lauret:EinsteinSolvandNilsolitonsCordobaConf2007} for the non-compact setting .

\begin{question*} Do all homogeneous spaces admit Einstein metrics?
\end{question*}

\begin{question*} Can the homogeneous Einstein spaces be classified?
\end{question*}

\begin{question*} What special properties do homogeneous Einstein metrics have?
\end{question*}

In the homogeneous setting,  Equation \ref{eqn: Einstein definition} is nothing more than a collection of quadratic equations and Einstein metrics are simply the (positive definite) real solutions.  The simplicity of this statement is deceiving and sits  in juxtaposition to the mountain of serious work invested over the past 70 years aimed at understanding these natural model spaces.

As is to be expected, the first efforts at addressing the existence question start with simply asking about constraints, both topological and Lie theoretic, of signed Ricci curvature.

Consider a homogeneous space $G/H$.  If $G/H$ admits a metric of positive Ricci curvature, by homogeneity Ricci curvature is bounded below by a positive constant, and so Myer's theorem yields that $G/H$ is compact.  Consequently, $G$ must be a compact, semi-simple group \cite{Jensen:HomogEinsteinSpacesofDim4}.  In Section \ref{sec: compact homog Einstein}, we discuss the current state of knowledge for the compact setting.

In the case of zero Ricci curvature, it turns out there is little variety.  In this case, the Riemannian metric must be flat and $G/H$ is isometric to $\mathbb R^k \times T^{n-k}$ \cite{AlekseevskiiKimelfeld:StructureOfHomogRiemSpacesWithZeroRicciCurv}.

When $G/H$ admits a metric of negative Ricci curvature, we know that $G$ and $G/H$ must be non-compact \cite{Bochner:VectorFieldsAndRicciCurvature}. This case has 
attracted substantial attention over the past 20 years.   We discuss this setting in Sections \ref{sec: non-compact homog Einstein}, \ref{sec: non-compact}, and \ref{sec: classification questions}.

There are two main, driving questions in the study of homogeneous Einstein spaces at present.  In the compact setting, Einstein metrics are not unique on a given homogeneous space.  On a given manifold, there are not even a finite number of homogeneous Einstein metrics.  However, if we fix the homogeneous presentation, then we have the following fundamental question.

\begin{openprob*} On a given compact homogeneous space $G/H$, are there only a finite number of $G$-invariant Einstein metrics (up to isometry and scaling)?
\end{openprob*}

In the non-compact setting, the Alekseekskii Conjecture has been the primary driver of research on Einstein spaces.

\begin{AC} If $G/H$ is connected and admits an Einstein metric with negative scalar curvature, then it is diffeomorphic to $\mathbb R^n$.
\end{AC}

Since the late 90s, there has been a flurry of activity towards resolving this conjecture.  A solution has recently been put forward \cite{Bohm-Lafuente:NonCompactEinsteinManifoldsWithSymmetry}.  We give an overview of the progress towards this resolution below.

Finally, any discussion of Einstein metrics would be incomplete without also considering their closely related cousins, the Ricci solitons.  Recall, a Riemannian metric $g$ on $M$ is a Ricci soliton if
	$$ric_g = cg +  L_Xg,$$
for some $c\in \mathbb R$ and some smooth vector field $X$ on $M$, where $ L_X$ is the Lie derivative.  These are important for two reasons.  First, both kinds of spaces arise naturally as fixed points/generalized fixed points of the (normalized) Ricci flow, see Section \ref{sec: ricci flow}.  Second, and most importantly for this article, in the non-compact homogeneous setting Ricci solitons are intimately coupled to Einstein metrics, see Section \ref{sec: non-compact homog Einstein}.

\subsection*{A note from the surveyor}
It is a privilege and a honor to have been asked  to write a survey on the topic of Einstein metrics on homogeneous spaces.  To do so within 20 pages is quite a challenge, but it is also an opportunity  for presenting a condensed narrative that should, hopefully, not overwhelm a newcomer to the field, giving them just enough to whet their appetite along with a sufficient amount of references  to help them dig deeper into wherever their interests lie and the whims of the moment point them.

Certainly, there will be many important contributions that have not be mentioned, some of these were powerful motivators whose technical results were eclipsed by other later works - still, they were necessary milestones along the way.  Other omissions will be due to this author's accidental omission or lack of awareness and knowledge on this author's part - for these omissions, I apologize in advance.

Finally, I have tried to fill the narrative with open questions that still need resolving.  Despite the significant progress made in recent years, these questions should demonstrate that the story is not fully written for homogeneous Einstein manifolds.

\subsection*{Acknowledgments} Special thanks go to Jorge Lauret for comments on an early draft of this manuscript and to the National Science Foundation who have supported this work under grant DMS-1906351.

\section{Compact Homogeneous Einstein spaces.}\label{sec: compact homog Einstein}

An excellent and robust survey of compact, homogeneous Einstein spaces can be found in \cite{Wang:EinsteinMetricsFromSymmetryAndBundleConstractionsASurvey}.  There the interested reader will find numerous explicit examples and a full treatment of the compact setting up till that point.  Since the writing of that survey, there have been advances in the form obtaining new and different kinds of examples of Einstein metrics on compact, homogeneous spaces, however,  the finiteness conjecture remains an open problem, see below.  Our treatment of the compact setting, in this section, should be viewed as an attempt to stimulate and motivate the reader to go deeper by exploring the survey by Wang.

The main tool for approaching Einstein metrics on compact, homogeneous spaces is the scalar curvature function restricted to the set of volume one metrics.  Let $\mathcal M_G$ denote the volume one $G$-invariant metrics on $G/H$ and 
	$$S: \mathcal M_G \to \mathbb R,$$
the  scalar curvature function.  Einstein metrics are critical points of this function \cite{Jensen:TheScalarCurvatureOfLeftInvariantRiemannianMetrics,
Nikonorov:TheScalarCurvFuctionalAndHomogEinsteinMetricsOnLieGroups}.

\begin{openprob}\label{openprob: finite number of Einstein compact homog} On a given homogeneous space $G/H$, are there a finite number of $G$-invariant Einstein metrics (up to isometry and scaling)?
\end{openprob}

Recall, we must restrict to a fixed  homogeneous structure $G/H$ as there do exist examples of infinitely many, inequivalent homogeneous structures on the same underlying manifold which admit $G$-invariant Einstein metrics \cite{Wang-Ziller:EinsteinMetricsOnPrincipalTorusBundles}.  For example, $S^2 \times S^3$ can be presented as $(SU(2)\times SU(2))/U(1)_{pq}$ with $p$ and $q$ relatively prime.

\begin{thm} \cite{Bohm-Wang-Ziller:AVariationalApproachforCompactHomogEinsteinMflds}  The moduli space of Einstein metrics on $\mathcal M_G$ has finitely many components, each of which is compact.
\end{thm}

Part of the challenge of the above is that the set of solutions to the Einstein equation $ric(g) = cg$ is not necessarily discrete.  For example, on a compact, semi-simple Lie group, one can pull-back any solution by an automorphism to have an isometric metric, which is of course also going to be Einstein.  On homogeneous spaces more generally, one has to contend with $N_G(H)/H$ which acts on $\mathcal M_G$.  (Note, we are not saying that all isometries arise this way.  That is a different, and interesting problem, in and of itself.  See, e.g., \cite{Gordon:RiemannianIsometryGroupsContainingTransitiveReductiveSubgroups,
GordonWilson:TheFineStructureOfTransitiveRiemannianIsometryGroups,
GordonWilson:IsomGrpsOfRiemSolv}.)  As a special case, the following finiteness conjecture has been made.

\begin{openprob}\cite{Bohm-Wang-Ziller:AVariationalApproachforCompactHomogEinsteinMflds} If G/H is a compact homogeneous space whose isotropy representation consists of pairwise inequivalent irreducible summands, e.g. when rank G = rank H, then the algebraic Einstein equations have only finitely many real solutions.
\end{openprob}

The above results and questions represent the leading edge of knowledge on the general compact setting.  Below, we give a few examples and sample both existence and non-existence results for particular homogeneous spaces.

\subsection{Examples of compact homogeneous Einstein spaces}

For the moment, let $G$ to be a Lie group, not necessarily compact, and consider $G/H$ where $H$ is compact.  Take an $Ad_H$-invariant decomposition $\mathfrak g = \mathfrak h \oplus \mathfrak p$.  If one prefers, one can choose $\mathfrak p$ to be the orthogonal complement to $\mathfrak h$ under the Killing form of $\mathfrak g$.  Then the set of $G$-invariant metrics on $G/H$ is one-to-one correspondance with the $Ad_H$-invariant inner products on $\mathfrak p \simeq T_{eH}G/H$.

The isotropy action of $H$ on $T_{eH}G/H$ is equivalent to the $Ad_H$-representation on $\mathfrak p$ and the simplest examples of homogenous Einstein spaces occur when this presentation is irreducible; that is, our homogeneous space is so-called isotropy irreducible.  The isotropy irreducible spaces are classified, see  \cite{Manturov:HomogeneousRiemannianSpacesWwithAnIrreducibleRotationGroup,
Wolf:TheGeometryAndStructureOfIsotropyIrreducibleHomogeneousSpaces}. These spaces include the following examples.

\begin{example}  Irreducible symmetric spaces are Einstein manifolds.
\end{example}

In the positive scalar curvature setting (i.e. compact), this includes the following.

\begin{example} Let $G$ be a compact, simple Lie group.  If $B$ is the Killing form of $\mathfrak g = Lie~G$, then we may consider the left-invariant $g_B$ which is $-B$ on $\mathfrak g \simeq T_eG$.
\end{example}

If $G$ is compact, semi-simple then we can build an Einstein metric by scaling the above on each factor. Likewise, for symmetric spaces, or a product of isotropy irreducible spaces, whose factors are of the same type (i.e. all of the same scalar curvature sign), we can scale the factors appropriately to put an Einstein metric on such a product.

\begin{thm}\cite{Jensen:TheScalarCurvatureOfLeftInvariantRiemannianMetrics,
DatriZiller:NaturallyRedMetricsAndEinsteinMetricsOnCompactLieGroups}
If $G$ is a compact semi-simple group other than $SU(2)$, then $G$ admits at least two Einstein metrics which are not equivalent up to scaling and isometry. The group $SU(2)$ admits a single Einstein metric.  
\end{thm}

In the more general, homogeneous setting, we have full information up to dimension 7, including an affirmation answer to Open Problem \ref{openprob: finite number of Einstein compact homog},  and partial information up to dimension 11.

\begin{thm}  The compact, homogeneous Einstein spaces are classified in dimension 7 and less.  For each homogeneous space $G/H$, there are a finite number of $G$-invariant Einstein metrics (up to scaling and isometry) on each $G/H$.
\end{thm}

Details of this classification can be found in \cite{Nikonorv:CompactHomogeneousEinstein7manifolds} which builds on work in dimesions 6 and less by \cite{Alekseevsky-Dotti-Ferraris:HomogeneosRicciPositive5manifolds,Nikonorov-Rodionov:CompactHomogeneousEinstein6manifolds}.  And, although there is no classification for higher dimensions, yet, one does know that Einstein metrics always exist up to dimension 11.

\begin{thm} \cite{Bohm-Kerr:LowDimensionalHomogeneousEinsteinManifolds}  Let $G$ be compact semi-simple with closed subgroup $H$.  If $G/H$ has dimension less than or equal to 11, then it admits a $G$-invariant Einstein metric.
\end{thm}

Among the examples given in \cite{Wang-Ziller:ExistenceAndNonexistenceOfHomogEinstein}  which do not admit Einstein metrics, there is an example of dimension 12, namely, $SU(4)/Sp(2)$, and so the above theorem is sharp.

Another case of note is homogeneous spaces $G/H$ for which the isotropy representation has only two summands.  The collection of such spaces admitting Einstein metrics has been classified \cite{Dickinson-Kerr:TheGeometryOfCompactHomogeneousSpacesWithTwoIsotropySummands,
He:CohomogeneityOneManifoldsWithASmallFamilyOfInvariantMetrics}, but the classification of Einstein metrics on these spaces remains open.

\subsection{General existence and non-existence results in the compact setting}
As Einstein metrics are critical points of the scalar curvature function (restricted to the volume one metrics), one might naturally look for maxima and minima of this function.

\begin{thm} \cite{Wang-Ziller:ExistenceAndNonexistenceOfHomogEinstein} If $H$ is a maximal proper subgroup of $G$, then $G/H$ admits a $G$-invariant Einstein metric.  Here the Einstein metric is a maximum of the scalar curvature function on $\mathcal M_1^G$.
\end{thm}

The work \cite{Wang-Ziller:ExistenceAndNonexistenceOfHomogEinstein} contains other criteria  for both the existence and non-existence of compact, homogeneous Einstein metrics.  Due to the above theorem, one must investigate how intermediate subgroups fit together, i.e. $K\subset H \subset G$,  and the relationship to non-extreme critical points of $S(g)$.  

One approach for guaranteeing or precluding the existence of non-extreme critical points is to build a simplicial complex out of intermediate algebras of $\mathfrak g$ with edges determined by inclusion of these algebras, and so on for building the faces.  Here the topology of this simplicial complex informs the higher energy critical values of $S$.  This is the approach taken in \cite{Bohm:HomogEinsteinMetricsAndSimplicialComplexes,Bohm:NonexistenceOfHomogeneousEinsteinMetrics}; in each case (i.e. guaranteeing the existence or non-existence of Einstein metrics), infinitely many examples are obtained.   Remarkably, one can use a finite amount of Lie algebra data to analyze different homogeneous spaces simultaneously.

\subsection{Results on compact, homogeneous Einstein manifolds since 2012}
In the ten years since Wang's extensive survey \cite{Wang:EinsteinMetricsFromSymmetryAndBundleConstractionsASurvey}, the bulk of new results for compact, homogeneous Einstein manifolds has been the construction of new types of Einstein metrics.  In some cases, these have been non-naturally reductive Einstein metrics  \cite{Yan-Deng:EinsteinMetricsOnCompactSimpleLieGroupsAttachedToStandardTriples,Chen-Chen-Deng:NonnaturallyReductiveEinsteinMetricsOnSOn,Yan-Deng:InvariantEinsteinMetricsOnCertainCompactSemisimple}, the classification of Einstein metrics for some special homogeneous spaces 
\cite{
Nikonorov:ClassificationOfGeneralizedWallachSpaces,
Chen-Nikonorov:InvariantEinsteinMetricsOnGeneralizedWallachSpaces}, and the generation of new examples of Einstein metrics on particular homogeneous spaces, sometimes with the aid of computer algebra systems, \cite{Arvanitoyeorgos-Sakane-Statha:HomogeneousEMetricsOnStiefelManifoldsAssociatedToFlagManifoldsWithTwoIsotropySummands,
Arvanitoyeorgos-Sakane-Statha:InvariantEinsteinMetricsOnSUnAndComplexStiefelManifolds,
Chrysikos-Sakane:HomogeneousEinsteinMetricsOnNonKahlerCSpaces}.  For a detailed treatment of Einstein metrics on flag manifolds, see \cite{Arvanitoyeorgos:ProgressOnHomogeneousEinsteinManifoldsAndSomeOpenProblems}.

\section{Non-compact, homogeneous Einstein spaces: solvmanifolds}\label{sec: non-compact homog Einstein}

As in the compact case, there are non-compact Lie groups which do admit Einstein metrics, namely solvable Lie groups.  If a solvable Lie group is unimodular, then it is Einstein if and only if it is flat \cite{Dotti:RicciCurvUnimodularSolv}.  Thus, a solvable Lie group admitting an Einstein metric of negative scalar curvature must be non-unimodular.  To date, all known examples of non-compact, homogeneous Einstein metrics on Lie groups occur on solvable groups.  As has recently been announced, this collection exhausts the set of non-compact, homogeneous Einstein spaces  \cite{Bohm-Lafuente:NonCompactEinsteinManifoldsWithSymmetry} and so we begin by focusing our attention on these groups.

As we will see in the sequel, the nilradicals of Einstein solvmanifolds have their own special properties.  For now, we call a nilpotent Lie group an Einstein nilradical if it is the nilradical of a solvable Lie group admitting an Einstein metric.  We call a nilpotent Lie algebra an Einstein nilradical if it's simply-connected Lie group is an Einstein nilradical.  (Simple-connectivity turns out to be necessary.)

\subsection{Examples of non-compact, homogeneous Einstein spaces}
As in the compact setting, symmetric spaces provide a familiar family of examples.  If $M$ is a non-compact, irreducible, symmetric space which is not flat, then $M= G/K$ where $G = Isom(M)$ is simple, $K$ is a maximal compact subgroup.  Even further, we have an Iwasawa decomposition $G = KAN$ with $S = AN$ being a solvable, simply-transitive group of isometries on $M$.  

The group $S$ above features several properties which, as we will see shortly, are the norm for solvable groups admitting Einstein metrics.
	\begin{enumerate}
	\item $S$ is simply-connected.
	\item The nilradical of $S$ is $N$ and $A$ is a complementary, abelian subgroup, i.e. $A\cap N = \{e \}$.
	\item The Lie algebra $\mathfrak s = Lie~S$ contains a special element $H\in\mathfrak a = Lie~A \subset \mathfrak s$ such that
		$$tr( ad_H \circ D) = tr~D \mbox{ for all derivations } D\in Der(\mathfrak s) .$$
	\end{enumerate}
For more details on the structure of Einstein solvmanifolds, see Section \ref{sec: structure Einstein solvmanifolds}.

By left-invariance, to build examples and study metrics on Lie groups, it suffices to work with Lie algebras and inner products on the Lie algebra.  We adopt this approach going forward.

\begin{example}  Let $\mathfrak n = \mathbb R^{n-1}$ be an abelian Lie algebra.  Consider the derivation $D$ which is the identity on $\mathfrak n$.  Then $\mathfrak s = \mathbb R D \ltimes \mathfrak n$ is the Lie algebra of a solvable group $S$ which admits an Einstein metric.  This produces the rank 1 symmetric space which is hyperbolic n-space $\mathbb H^n$.
\end{example}

\begin{example} Let $\mathfrak n$ be the Heisenberg Lie algebra algebra.  Let $\mathfrak z$ be the center of $\mathfrak n$ and write $\mathfrak n = \mathfrak z \oplus \mathfrak v$ for some vector complement $\mathfrak v$ of $\mathfrak z$.  Define a derivation $D$ of $\mathfrak n$ which is the identity on $\mathfrak v$ and twice the identity on $\mathfrak z$.  Then $\mathfrak s = \mathbb R D \ltimes \mathfrak n$ admits an inner product such that the corresponding metric on $S$ is Einstein.  This is the rank 1 symmetric space  complex hyperbolic space
\end{example}

In the above two examples, using any other derivation will product a group that cannot admit an Einstein metric.  Unlike in the compact setting, not all (non-unimodular) solvable groups admit left-invariant Einstein metrics in low dimensions due to the necessary, algebraic structures above.  However, in low dimensions, those are the only constraints.

\begin{thm} \cite{Lauret:FindingEinsteinSolvmanifoldsByAVariantialMethod,Will:Rank1EinsteinSolvOfDim7} Up to dimension 6, every  nilpotent Lie algebra is an Einstein nilradical.
\end{thm}

In dimension 7, there are examples of non-Einstein nilradicals, but also continuous families of Einstein nilradicals; these are classified in dimension 7 
\cite{LauretWill:EinsteinSolvExistandNonexist,
Fernandez-Culma:ClassificationOfNilsolitonMetricsInDimensionSeven}.  At this point, we start to see more divergence in the theory for compact and non-compact spaces.  In the compact setting, in any given dimension, the number of compact, semi-simple Lie groups is finite; in the solvable setting, beyond low dimensions you have large dimension continuous families of solvable groups.

Another difference between the compact and non-compact settings is that for non-unimodular solvable Lie groups, Einstein metrics are not critical points of the scalar curvature function.  However, one can realize them as critical points  of a modified scalar curvature function \cite{Heber,Lauret:StandardEinsteinSolvAsCriticalPoints}.  

As we will see in Section \ref{sec: structure Einstein solvmanifolds}, Einstein metrics on solvable Lie groups are unique up to scaling and isometry, yet another divergence from the compact realm.

\begin{thm}\cite{Eberlein:prescribedRicciTensor,Jablo:Thesis}  A generic 2-step nilpotent Lie algebra is an Einstein nilradical.
\end{thm}

The first examples of continuous families of non-Einstein nilradicals appeared  in \cite{Will:CurveOfNonEinsteinNilradicals}.  This was followed by other techniques for building (continuous) families of both Einstein and non-Einstein nilradicals  \cite{Payne:ExistenceofSolitononNil,Arroyo:FiliformNilsolitonsOfDimension8,
Jablo:ModuliOfEinsteinAndNoneinstein,Payne:GeometricInvariantsForNilpotentMetricLieAlgebrasWithApplicationsToModuliSpacesOfNilsolitonMetrics}.  See \cite{Lauret-Oscari:OnNonsingular2stepNilpotentLieAlgebras,
Oscari:OnTheExistenceOfNilsolitonsOn2stepNilpotentLieGroups} for a deeper investigation of non-Einstein nilradicals among 2-step nilpotent algebra. 
In the presence of additional curvature conditions, one can classify the Einstein solvmanifolds in low dimensions.  For example, if one requires negative sectional curvature, then these are classified in dimensions up to 7 \cite{NikitenkoNikonorov:Six-dimensionalEinsteinSolvmanifolds}\cite{Nikitenko:Seven-dimensionalHomogeneousEinsteinManifoldsOfNegativeSectionalCurvature}.

Unlike in the compact setting, we have a convenient reduction to irreducible algebras.

\begin{thm}\cite{Jablo:DetectingOrbitsAlongSubvarietiesViaTheMomentMap,Nikolayevsky:EinsteinSolvmanifoldsandPreEinsteinDerivation} Let $\mathfrak n = \mathfrak n_1 \oplus \mathfrak n_2$ be a nilpotent Lie algebra which is a direct sum of ideals.  If $\mathfrak n$ is an Einstein nilradical, so are the factors $\mathfrak n_1$ and $\mathfrak n_2$.
\end{thm}

Of note, it is not assumed apriori that the Einstein metric makes $\mathfrak n_1$ and $\mathfrak n_2$ orthogonal.  Instead, this is true after the fact.

\subsection{Structure results for Einstein solvmanifolds}\label{sec: structure Einstein solvmanifolds}

Presently, there are few techniques for precluding the existence of Einstein metrics.  First, among Einstein metrics, zero scalar curvature ($c=0$) corresponds to unimodular Lie groups and negative scalar curvature ($c<0$) corresponds to non-unimodular Lie groups \cite{Dotti:RicciCurvUnimodularSolv}.  As the Ricci flat homogeneous spaces are precisely the flat spaces \cite{AlekseevskiiKimelfeld:StructureOfHomogRiemSpacesWithZeroRicciCurv}, we focus our attention on the setting of negative scalar curvature.

To demonstrate the stark contrast between the compact and non-compact settings, we begin with the following.

\begin{thm} \cite{Heber}\label{thm: Heber on uniqueness of Einstein} A left-invariant (standard) Einstein metric on a non-unimodular solvable Lie group is unique up to scaling and isometry.
\end{thm}

A solvable Lie group with left-invariant metric is called standard if the corresponding Lie algebra $\mathfrak s$ with inner product satisfies $\mathfrak s = \mathfrak a \oplus \mathfrak n$ where $\mathfrak n$ is the nilradical and $\mathfrak a = \mathfrak n ^\perp$ is an abelian Lie algebra.  The standard condition turns out to always be satisfied (see below) and so we place it in parentheses.

\begin{remark} In the theorem above, we actually know that any two Einstein metrics are equivalent up to scaling and pull-back by an automorphism.  While pull-back by an automorphism will always give isometric metrics on a Lie group, in general not all isometries arise this way for all Lie groups, even among solvable Lie groups.  See \cite{GordonWilson:IsomGrpsOfRiemSolv}.
\end{remark}

We have a recipe for building Einstein solvmanifolds which goes as follows.  Take as ingredients:
	\begin{enumerate}
	\item A nilpotent Lie algebra $\mathfrak n$ with a so-called Ricci soliton metric (see Equation \ref{eqn: Ricci soliton defintion}).
	\item An abelian subalgebra $\mathfrak a \subset Der(\mathfrak n)$ of fully-reducible operators whose real parts are non-trivial.
	\item Require the existence of an element $H\in\mathfrak a$ such that $D=ad~H \in Der(\mathfrak n)$ satisfies
		\begin{equation}\label{Eqn: pre-Einstein derivation}
		tr(D\phi) = tr(\phi) \mbox{ for all } \phi\in Der(\mathfrak n)
		\end{equation}
	with positive eigenvalues.
	\end{enumerate}
We can assume, by conjugating $\mathfrak a$ by an automorphism of $\mathfrak n$, that $\mathfrak a$ consists of normal operators relative to the Ricci soliton inner product on $\mathfrak n$.  On $\mathfrak a$, use the metric $\ip{A,B} = tr (S(ad_A) \circ S(ad_B))$, where $S(X) = \frac 1 2 (X+X^t)$ denotes the symmetric part of the normal operator $X$.  On $\mathfrak n$ use the Ricci soliton inner product.  Extend these to an inner product on  the Lie algebra $\mathfrak s = \mathfrak a \ltimes \mathfrak n$ so that $\mathfrak a \perp \mathfrak n$.  We now have the metric Lie algebra of a solvable Lie group with left-invariant Einstein metric.  Remarkably, this construction exhausts the class of Einstein solvmanifolds.

\begin{thm} \cite{Lauret:EinsteinSolvmanifoldsAreStandard,Lauret:SolSolitons} \label{thm: Lauret Einstein solv are standard}
All Einstein solvmanifolds are standard and arise via the construction above.

\end{thm}

We note that partial progress on the problem of showing Einstein solvmanifolds are standard was made by others, see e.g. \cite{Nikonorov:OnEinsteinExtensionsOfNilpotentMetricLieAlgebras} and other references of \cite{Lauret:EinsteinSolvmanifoldsAreStandard}.  This theorem is  the root of all non-existence results for Einstein metrics on solvable Lie groups.  We note a few important consequences.

First, the special derivation $D$ defined in Equation \ref{Eqn: pre-Einstein derivation} must have eigenvalues which are intergers (up to scaling) and so $\mathfrak n$ must be $\mathbb N$-graded.  Recall, a (nilpotent) Lie algebra is $\mathbb N$-graded if
	$$\mathfrak n  = \mathfrak n_1 \oplus \dots \oplus \mathfrak n_k \mbox{ with } [\mathfrak n_i,\mathfrak n_j] \subset \mathfrak n_{i+j}.$$
Here some of the $\mathfrak n_i$ may be trivial.  This is equivalent to the existence of a derivation with eigenvalues which are positive integers. 

\begin{example}  Characteristically nilpotent Lie algebras cannot be Einstein nilradicals.
\end{example}

Recall, a nilpotent Lie algebra is called characteristically nilpotent if $Der(\mathfrak n)$ is nilpotent.  
Going beyond the necessity of the existence of a positive derivation, one must have a derivation satisfying Equation \ref{Eqn: pre-Einstein derivation} and this has served as  inspiration for the pre-Einstein derivation as defined in \cite{Nikolayevsky:EinsteinSolvmanifoldsandPreEinsteinDerivation}.  This pre-Einstein derivation turns out to play a deeper role in the geometry of solvmanifolds, is unique up to conjugation in $Der(\mathfrak n)$, and is a necessary ingredient in understanding maximal symmetry, see Section \ref{sec: maximal symmetry}.

\begin{thm}\cite{Nikolayevsky:EinsteinSolvmanifoldsandPreEinsteinDerivation} If $\mathfrak n$ is an Einstein nilradical, then its pre-Einstein derivation must be positive, i.e. have positive eigenvalues.
\end{thm}

\begin{openprob}\label{openprob: is there list of algebraic invariants for Einstein condition}  Is there a list of algebraic invariants which completely determines when a solvable algebra admits an Einstein metric?
\end{openprob}

Despite not having  a list of algebraic invariants, we do know that the existence of an Einstein metric on a solvable Lie group is intrinsic to the Lie algebra and is a `local problem' on the space of left-invariant metrics.  We are careful to point out that finding the Einstein metric itself is a global problem on the space of metrics.  In \cite{Jablo:ConceringExistenceOfEinstein} it was shown that one can determine the existence of an Einstein metric by measuring 1) algebraic invariants and 2) local deformations of any choice of left-invariant metric.  Like most of the results on existence of Einstein metrics, this results uses the robust tool of Geometric Invariant Theory described in Section \ref{subsec: GIT}

\section{Non-compact, homogeneous Einstein spaces: tools and structure results}\label{sec: non-compact}
A central tool in the study of non-compact, homogeneous Einstein spaces of negative scalar curvature is Geometric Invariant Theory.  Its introduction to homogeneous Riemannian geometry was in the work \cite{Heber}.  We motivate its introduction with a close look at the Ricci tensor of a Lie group with left-invariant metric.  Note, everything that follows can and has been extended to homogeneous spaces.

On a Lie group with left-invariant metric, one can consider the Ricci tensor and decompose it as follows
	$$Ric = M -\frac 1 2 B - S(ad_H),$$
where $B(X,X)= - tr~(ad_X)^2$ is the Killing form, $H$ a mean curvature vector satisfying $\ip{H,X}=tr~ad_X$ for all $X\in\mathfrak g$, with $S(ad_H)$ the symmetric part of $ad_H$, and $M$ a mysterious component \cite{Besse:EinsteinMflds}.  
This mysterious tensor $M$ is the moment map for the natural, change of basis action on the space of Lie brackets; to our knowledge, this crucial observation was first made by Lauret \cite{Lauret:DegenerationsOfLieAlgebrasAndGeometryOfLieGroups,
Lauret:MomentMapVarietyLieAlgebras}.

\subsection{Geometric Invariant Theory as a tool}\label{subsec: GIT}
Here we describe a tool which has been exploited in the nilpotent, solvable, and general homogeneous setting for understanding Einstein and Ricci soliton geometries.  We describe everything in setting of Lie groups with left-invariant metrics for simplicity, however the reader should note that the ideas all extend to the full homogeneous setting  \cite{LauretLafuente:StructureOfHomogeneousRicciSolitonsAndTheAlekseevskiiConjecture}.

Consider a Lie group $G$ with left-invariant metric $g$.  If $G$ is simply-connected, all the data for $(G,g)$ is encoded in the metric Lie algebra $(\mathfrak g, \innerprod)$. Here the inner product on $\mathfrak g \simeq T_eG$ is the restriction of $g$ to $ T_eG$.  A metric Lie algebra is three pieces of data:  an underlying vector space $\mathbb R^n$, a Lie bracket $\mu = [\cdot , \cdot ]$, and an inner product $\innerprod$.  We will let $G_{\mu,\innerprod}$ represent the simply-connected Lie group with metric Lie algebra $\{\mathbb R^n , \mu , \innerprod \}$.

On the one hand, we can vary the inner product to obtain all possible left-invariant geometries on $G$ via
	$$\phi\cdot \innerprod = \ip{\phi \cdot , \phi \cdot } \mbox{ for } \phi \in GL(n,\mathbb R);$$
here $GL(n,\mathbb R)$ is acting on the space of inner products which is the open set $\mathcal P$ of positive definite, symmetric bilinear forms.
On the other hand, we can vary the Lie structure via the change of basis action  as follows
	$$\phi\cdot \mu = \phi \mu (\phi^{-1} \cdot, \phi^{-1}\cdot ) \mbox{ for } \phi \in GL(n,\mathbb R);$$
here $GL(n,\mathbb R)$ is acting on the space of $\mathbb R^n$-valued anti-symmetric, bilinear forms $\wedge ^2 (\mathbb R^n)^*\otimes \mathbb R^n$.  The set of Lie brackets is an algebraic variety in this vector space as the Jacobi condition is polynomial.  These two perspectives are equivalent and we have that the following Lie groups with left-invariant metrics are isometric
	$$G_{\mu, \phi\cdot\innerprod} \simeq G_{\phi\cdot\mu, \innerprod} .   $$
The isometry between these Lie groups with left-invariant metrics arises from lifting the Lie algebra isomorphism $\phi: (\mathbb R^n, \mu) \to (\mathbb R^n, \phi\cdot \mu)$ which is simultaneously an isometry of the vector spaces with inner products.

The slight shift in perspective to varying the Lie bracket turns out to be quite powerful.  We demonstrate the usefulness of this shift in perspective for nilpotent Lie groups.  For a simply-connected, nilpotent Lie group $N$, the set of isometry classes in $\mathcal P$ are precise the orbits of $Aut(N)$ \cite{GordonWilson:IsomGrpsOfRiemSolv}.  However, in the space of Lie brackets the isometry classes become the orbit of $O(n)$.

Additionally, by working in the space of Lie brackets, we have a natural way of compactifying the set of left-invariant metrics.  For a nilpotent Lie bracket $\mu$, $sc(N_\mu) = -\frac{1}{4}|\mu|$ and so normalizing to left-invariant metrics of $sc = -1$ amounts to restricting oneself to a sphere in $\wedge^2(\mathbb R^n)^*\otimes \mathbb R^n$.  The compactification of the set of left-invariant metrics on a nilpotent Lie group is then 
	$$ \overline{GL(n,\mathbb R)\cdot \mu} \cap S,$$
where $S$ is a sphere of radius 4 in $\wedge^2(\mathbb R^n)^*\otimes \mathbb R^n$.  Obviously, in the boundary, we have non-isomorphic Lie groups that appear in the compactification.  This turns out to be an important, and useful, fact.

Further, for nilpotent Lie groups $G_{\mu,\innerprod}$, we have 
	$$Ric = M_\mu$$
where $Ric$ is the $(1,1)$-Ricci tensor and $M_\mu$ is the moment map of the $GL(n,\mathbb R)$ action on $\wedge^2(\mathbb R^n)^*\otimes \mathbb R^n$.  This follows as nilpotent groups being unimodular implies $H=0$  and nilpotentcy implies the Killing form $B$ vanishes.  

From Geometric Invariant Theory \cite{Kirwan}, there is a natural stratification of $\wedge^2(\mathbb R^n)^*\otimes \mathbb R^n$ which plays a crucial role in the possible geometries a Lie group can admit.  These stratifications for non-nilpotent groups play an essential role in the structure results for general, non-compact, homogeneous Einstein manifolds \cite{LauretLafuente:StructureOfHomogeneousRicciSolitonsAndTheAlekseevskiiConjecture}.  We layout the basics of the GIT stratification to motivate the interested reader to dig deeper.

For the $GL(n,\mathbb R)$ action on $V = \wedge^2(\mathbb R^n)^*\otimes \mathbb R^n$, we have an induced action of $\mathfrak{gl}(n,\mathbb R)$ on $V$ which will be denoted by $\pi$.  We implicitly define the function 
	$$ m : V = \wedge^2(\mathbb R^n)^*\otimes \mathbb R^n \to \mathfrak{gl}(n,\mathbb R)$$
via 
	$$\ip{ m(\mu) , X } = \frac{1}{|\mu|^2}\ip {\pi(X) \mu,\mu}.$$
Note, $m(\mu) = \frac{1}{|\mu|^2}M_\mu$.  The function $F(\mu) = | m(\mu)|^2$ has a finite number of critical values.  Let $\mathcal C$ denote the critical points of $F$.  There is a finite collection of diagonal metrics $\mathcal B$ such that if $\mu$ is a critical point, then $m(\mu) = k\beta k^{-1}$ for some $\beta \in \mathcal B$ and some $k\in O(n)$.  

Let $\mathcal C_\beta$ be the critical points with $m(\mu)$ conjugate to $\beta$.  We define $\mathcal S_\beta$ to be the points in $V$ which flow under the negative gradient flow to $\mathcal C_\beta$.  In this way, the critical points are precisely the minima of $F$ on the stratum $\mathcal S_\beta$.

As to be expected, for the particular representation of $GL(n,\mathbb R)$ above, these strata enjoy extra properties since the points in our space correspond to algebras.  For example, we have
	\begin{enumerate}
	\item For $\mu \in \mathcal S_\beta$, $\ip{[\beta,D],D}\geq 0$ for all $D\in Der(\mu)$.
	\item If $\mu$ is nilpotent, then the critical points of $F$ are precisely the Ricci soliton metrics, cf. Equation \ref{eqn: Ricci soliton defintion}.
	\item If $\mu$ is solvable, then $GL(n,\mathbb R)\cdot \mu \cap \mathcal C \not = \emptyset$ if and only if $G_\mu$ admits a left-invariant Einstein metric.  Here the critical point is not the Einstein metric, but one can obtain the Einstein metric by simply dilating the orthogonal complement to the nilradical.
	\end{enumerate}
The GIT stratification is a corner stone to obtaining the general  structure results for homogeneous, Einstein spaces.

\subsection{Structure results for non-compact, homogeneous Einstein spaces in general}
We now consider the full setting of non-compact, homogeneous Einstein spaces and not just the solvmanifold setting.  

\begin{thm} Suppose $G$ is non-compact and unimodular with compact subgroup $H$.  If $G/H$ admits a $G$-invariant Einstein metric, then it is Ricci flat and hence flat.
\end{thm}

From this, we see that the lion's share of our interest will be for groups $G$ which are non-unimodular.  Using Geometric Invariant Theory in the full homogeneous setting and exploiting the GIT-stratification for our particular cases, we have the following.

\begin{thm} 
\cite{LauretLafuente:StructureOfHomogeneousRicciSolitonsAndTheAlekseevskiiConjecture,
JP:TowardsTheAlekseevskiiConjecture,
Arroyo-Lafuente:TheAlekseevskiiConjectureInLowDimensions}
Let $M$ be a non-compact, homogeneous Einstein space of negative scalar curvature.  Let $G$ denote the connected component of the identity of the isometry group.  Here $M=G/H$ with $H$ compact.  We have a Levi decomposition $G= G_1G_2$ where $G_1$ is a maximal semi-simple and $G_2$ the radical (i.e. maximal, normal solvable Lie subgroup) such that
	\begin{enumerate}
	\item $H = H_1H_2$ with $H_1 = H \cap G_1$ and $H_2=H\cap G_2$
	\item $G_2 = A_2 H_2 N_2$ where $N_2$ is the nilradical of $G$, $A_2H_2$ is an abelian subgroup with $A_2$ acting on $N_2$ with $ad~\mathfrak a_2$ consisting of reductive operators with real eigenvalues
	\item The group $A_2N_2$ with the submanifold geometry is an Einstein manifold itself
	\item writing $G_1 = G_c G_{nc}$ as a product of $G_c$ the product of simple compact, normal subgroups and $G_{nc}$ the product of simple non-compact, normal subgroups we have $G_c \subset H$.
	\end{enumerate}	
\end{thm}

Essentially, the above allows us to remove the compact, normal subgroups from $G_1$ and reduce $G_2$ to the setting of completely solvable groups.  At this point, what is left to resolve is how large is the isotropy $H_2 \subset G_1$.  If $H_2$ is a maximal compact subgroup, then one has that $G/H$ is a solvmanifold.

\subsection{Overview of the resolution of the conjecture}
The Alekseevskii Conjecture asserts that a (connected) non-compact, homogeneous Einstein space with negative scalar curvature is diffeomorphic to Euclidean space.  This conjecture was known to hold in low dimensions, namely dimensions 4, 5, and 7 from the works \cite{Jensen:HomogEinsteinSpacesofDim4,
Nikonorov:NoncompactHomogEinstein5manifolds,
Arroyo-Lafuente:TheAlekseevskiiConjectureInLowDimensions}, respectively.  Partial results in dimensions up to 10 were obtained in \cite{Arroyo-Lafuente:TheAlekseevskiiConjectureInLowDimensions,
Berichon:TheAlekseevskiiConjectureIn9and10dimensions}.

We now layout a strategy for resolving the conjecture that is ahistorical - i.e. instead of presenting the results as they appear chronologically, we present a narrative  which exploits the benefit that hindsight affords, with necessary structural results having appeared and been refined at various points in time.

First, we note that non-compact  Einstein solvmanifolds of negative scalar curvature (i.e. $c<0$) are simply-connected.  By solvmanifold, we mean that it has a transitive, solvable group of isometries.  A natural starting place would be to transform the original Alekseevskii Conjecture which is a topological statement into one whose outcome is Lie theoretic.

\begin{thm}\cite{Bohm-Lafuente:HomogeneousEinsteinMetricsOnEuclideanSpacesAreEinsteinSolvmanifolds} Supose $(G/H,g)$ is a simply-connected homogeneous Einstein space with negative scalar curvature.  If $G/H$ is diffeomorphic to $\mathbb R^n$, then $G/H$ is a solvmanifold.
\end{thm}

Next, one can reduce to the simply-connected setting.

\begin{thm}\cite{Jablo:StronglySolvable}  Let $(G/H,g)$ be a homogeneous Einstein space whose simply-connected cover is a solvmanifold.  Then $G/H$  is a solvmanifold and the quotient $\widetilde{G/H} \to G/H$ is trivial.
\end{thm}

It then remains to prove the Alekseevskii Conjecture among simply-connected manifolds. This is the assertion in the recent work \cite{Bohm-Lafuente:NonCompactEinsteinManifoldsWithSymmetry}.  A new suite of tools is developed there which we will not comment on in this survey and we refer the interested reader to that new work.

\section{Classification questions}\label{sec: classification questions}
In the compact setting, the classification of homogeneous Einstein metrics seems out of reach beyond low dimensions, see Section \ref{sec: compact homog Einstein}.  This is in part due to the fact that most compact homogeneous spaces admitting Einstein metrics admit more than one.   However, in the non-compact setting, i.e. the setting of solvmanifolds, Einstein metrics are unique (up to isometry and scaling) on a given solvmanifold when they exist, see Theorem \ref{thm: Heber on uniqueness of Einstein}.  And so  it is reasonable to ask if we can classify the solvable Lie groups/homogeneous spaces that admit an Einstein metric.  

In dimensions 7 and greater, the number of solvmanifolds admitting Einstein metrics is infinite \cite{LauretWill:EinsteinSolvExistandNonexist} and we have to ask what it even means to classify these spaces. Any reasonable classification should 
	\begin{enumerate}
	\item be a list of the objects without redundancy and
	\item provide a means for determining if an object you bring to the list is either on it or not.
	\end{enumerate}
Going further, we might even ask
	\begin{enumerate}
	\setcounter{enumi}{2}
	\item to identify exactly where in the list an object is.
	\end{enumerate}
The second point is the challenge as it is a difficult problem to know when any two given Lie algebras are isomorphic and the third point might seem impossible.  However, both can be achieved in reasonable sense for our problem.

\subsection{The classification list}
We follow \cite{Will:TheSpaceOfSolsolitonsInLowDimensions} to produce a list for the classification of Einstein solvmanifolds.

The first step is to classify the Einstein nilradicals or equivalently nilsolitons, cf. Theorem  \ref{thm: Lauret Einstein solv are standard}.  From Geometric Invariant Theory (see Section \ref{subsec: GIT}), we know that the set of nilpotent Lie groups with Ricci soliton metrics is precisely the set of critical points of the $F(\mu) = |m(\mu)|^2$, where $m(\mu)$ is the (normalized) moment map.  Denoting the critical points of $F$ by $\mathcal C$, we see that for  $\mu \in \mathcal C$, 
	$$GL(n,\mathbb R)\cdot \mu \cap \mathcal C = \mathbb R( O(n)\cdot \mu).$$
In this way, we see that nilsoliton metrics on a given Lie algebra are unique up to scaling and isometry, cf.  \cite{Lauret:RicciSolitonHomogeneousNilmanifolds}.  Thus, the set of Einstein nilradicals is in one-to-one correspondence with
	$$\mathcal C\cap \mathcal N/(\mathbb R\times O(n)),$$
where $\mathcal N$ is the set of nilpotent Lie brackets in $V = \wedge^2 (\mathbb R^n)^* \otimes \mathbb R^n$.  

At each metric nilpotent Lie algebra $\mu$ above, we can consider a maximal abelian, symmetric subalgebra $\mathfrak a$ of $Der(\mu)$.  This subalgebra is unique up to conjugation in $Aut(\mu)$ and so is essentially unique.  Suppose $\dim \mu = n$ and $\dim \mathfrak a = a$.  For each $1\leq r \leq a$, we may consider the space of solvsolitons of dimension $r+n$ with a given nilradical $\mu$.  This space is parameterized by a quotient of the Grassmannian
	$$Gr_r(\mathfrak a)/W$$
where $W$ is an explicit finite group \cite{Will:TheSpaceOfSolsolitonsInLowDimensions}.  In this way, one has a `list' of the solvsoliton metrics or equivalently the algebras admitting such metrics.  To obtain a list of the Einstein solvmanifolds, one simply restricts to the subspace in $Gr_r(\mathfrak a)$ which contain the pre-Einstein derivation, see Equation \ref{Eqn: pre-Einstein derivation}.  In a sense, this resolves the first point of having a classification list.

\subsection{Determining if you are on the list}
Our next job is to determine whether or not a Lie algebra is on our list when one is handed to us.  This can be achieved, in principle, but not in the most ideal of ways.  Ideally, one would have a list of algebraic invariants which resolves the question of when a given Lie algebra admits a soliton or Einstein metric.  This is been a long-standing open question for Einstein solvmanifolds \cite{LauretWill:EinsteinSolvExistandNonexist}.

At this point in time,  the best we have  is a combination of 
	\begin{enumerate}
	\item algebraic invariants/measurements and 
	\item `local measurement' using any initial metric.
	\end{enumerate}
This has been achieved \cite{Jablo:ConceringExistenceOfEinstein}. See the discussion after the Open Question \ref{openprob: is there list of algebraic invariants for Einstein condition} for details.

\subsection{Determining where on the list you are}
Since our classification list is actually a list of metric Lie algebras (with the metric of preference), satisfying the third point is a matter of being able to find the soliton/Einstein metric on a given Lie algebra.  Presently, there are few tools for finding these special metrics.  In the Einstein setting, Ricci flow is known to take any initial metric and produce the Einstein metric at time infinity, upon running the flow \cite{Bohm-Lafuente:TheRicciFlowOnSolvmanifoldsOfRealType}.

\begin{openprob}  Are there algebraic techniques for recovering the Einstein metric of a given solvable Lie group when one is known to exist?
\end{openprob}

\section{Special properties of Einstein spaces}
Given that Einstein metrics  on solvable Lie groups are unique (up to isometry and scaling), one might expect these metrics to have special properties.  This is certainly the case.

\subsection{Ricci flow}\label{sec: ricci flow}

Since, by definition, Einstein metrics satisfy $ric(g) = c g$, these metrics essentially do not change under the  Ricci flow - under the flow they simply dilate.  Recall, the Ricci flow is the differential equation on the space of Riemannian metrics given by
	$$\frac{\partial}{\partial t} g = -2 ric(g). $$
One can normalize by rescaling to the Ricci flow to fix  scalar curvature; for this normalized flow, Einstein metrics are then fixed points.  We note that for compact manifolds, the isometry group is preserved and so homogeneity is preserved under the Ricci flow.  For non-compact manifolds, the situation is more delicate \cite{Kotschwar:BackwardsUniquenessRF}, so one has to choose to restrict the Ricci flow to the set homogeneous metrics.  On the set of homogeneous metrics, the PDE becomes an  ODE and solutions to the Ricci flow are unique.

\begin{question} Are the fixed points of the normalized Ricci flow (i.e. Einstein metrics) stable?
\end{question}

Recall, for unimodular Lie groups $G$ (including compact Lie groups) and on the space of $G$-invariant metrics on $G/H$, the Ricci tensor is the gradient of the scalar curvature function.  Considering there are Einstein metrics of various energy levels for $G$ compact, we cannot have all Einstein metrics be stable.

\begin{thm}\cite{Kroncke:StabilityOfEinsteinMetricsUnderRicciFlow} $\mathbb CP^n$ with its Fubini-Study metric (which is Einstein) is not stable.
\end{thm}

In the work \cite{Kroncke:StabilityOfEinsteinMetricsUnderRicciFlow}, criteria are given for determining when a compact Einstein manifold is stable.   This work has been extended to other symmetric spaces \cite{Hall-Murphy-Waldron:CompactHermitianSymmetricSpaces}.  On a fixed homogeneous space $G/K$, and looking at only $G$-invariant metrics, significantly more is known; we direct the interested reader to the recent work \cite{Lauret:OnTheStabilityOfHomogeneousEinsteinManifolds}.

For non-compact spaces, we are able to have more hope as Einstein metrics are unique here.

\begin{thm}\cite{Bohm-Lafuente:TheRicciFlowOnSolvmanifoldsOfRealType} Let $S$ be a solvable Lie group with left-invariant Einstein metric $g$.  Among homogeneous metrics, $g$ is stable under the Ricci flow.  In fact, the Ricci flow starting at any left-invariant metric evolves to the Einstein metric.
\end{thm}

It should be no surprise that the proof of this result relies on the robust structure theory and GIT-stratification on the space of Lie brackets, cf. Section \ref{subsec: GIT}.  Since the Ricci flow preserves the isometry group of a metric, this gives a dynamical proof of maximal symmetry since the Einstein metric appears in the limit of the flow - see Theorem \ref{thm: Einstein solv have maximal symmetry}.

On homogeneous spaces, in general, the Ricci flow can be converted to a flow on the space of brackets \cite{Lauret:RicciFlowOfHomogeneousManifoldsAndItsSolitons}.  This approach allows one to prove results on the longterm existence of the flow and more.  For example, we have the following.

\begin{thm}
\cite{Lauret:RicciFlowForSimplyConnectedNilmanifolds,
Lauret:RicciFlowOfHomogeneousManifoldsAndItsSolitons} On nilpotent Lie groups with left-invariant metrics, the Ricci flow always evolves towards a Ricci soliton (defined below) - on a possibly different nilpotent group.
\end{thm}

The result above exploits GIT upon observing that  the Ricci flow is essentially the negative gradient flow of the norm-squared of the moment map on the space of brackets.  In low dimensions, we also refer the interested reader to \cite{GlickensteinPayne:RicciFlow3dimUnimodLieAlgs}.  More stability results are presented for nilpotent Lie groups with Ricci solitons below.

\subsection{Homogeneous Ricci solitons}
From the perspective of the Ricci flow, if one is interested in fixed points, one should also consider Ricci solitons as they are generalized fixed points in the sense that `up to diffeomorphism' they are fixed.  A metric $(M,g)$ is called a Ricci soliton if 
	\begin{equation}\label{eqn: Ricci soliton defintion}
	ric(g) = cg +  L_X
	\end{equation}		
where $L_X$ is the Lie derivative of  some smooth vector field $X\in \mathfrak X(M)$.  This is equivalent to $g$ being the initial point of a solution to the Ricci flow of the form $g(t) = c(t) \phi(t)^* g$ for a 1-parameter family of diffeomorphisms $\phi(t)$ and a 1-parameter family of constants $c(t) >0$.

\begin{question}  What are the homogeneous Ricci solitons?
\end{question}

Surprisingly, the answer to this question is quite simple and they are all born from the homogeneous Einstein metrics.

\begin{thm}\cite{Naber:NoncompactShrinking4SolitonsWithNonnegativeCurvature,Petersen-Wylie:OnGradientRicciSolitonsWithSymmetry}
 If $G/H$ admits a Ricci soliton with non-negative cosmological constant $c\geq 0$, then $G/H$ is (locally) the product of a compact homogeneous Einstein space with a Euclidean factor.
\end{thm}

Compact solitons are  necessarily Einstein with $c\geq 0$ \cite{Jablo:HomogeneousRicciSolitons} and so we  focus on the setting of non-compact, homogeneous solitons with $c < 0$.

\begin{thm}\cite{Jablo:HomogeneousRicciSolitonsAreAlgebraic} Consider a (non-compact) homogeneous Ricci soliton $M$ with negative cosmological constant $c$.  Let $G$ denote the isometry group of $M$, then our Ricci soliton is algebraic with respect to $G$, i.e. there is a symmetric derivation $D\in Der(\mathfrak g)$ such that the vector field $X$ in Equation \ref{eqn: Ricci soliton defintion} is generated by the 1-parameter family of automorphisms $exp(tD)$.  Even further, the $(1,1)$-Ricci tensor is of the form 
	$$Ric = cId + D_\mathfrak p,$$
where $\mathfrak g = \mathfrak k \oplus \mathfrak p$, $T_e G/K \simeq \mathfrak p$, and $D_\mathfrak p$ is  the projection of $D$ onto $\mathfrak p$.
\end{thm}

\begin{remark}
Of note, the map $D_\mathfrak p$ above is symmetric.  This subtle detail turns out to be essential for various structure results in the classification Einstein and soliton metrics on non-compact homogeneous spaces.
\end{remark}

As in the Einstein case, one can ask about the stability of homogeneous Ricci solitons under the Ricci flow.  In \cite{WilliamsWu:DynamicalStabilityOfAlgebraicRicciSolitons}, it is shown that linear stability of solitons implies dynamical stability.  In \cite{JPW:LinearStabilityOfAlgebraicRicciSolitons,JPW2:OnTheLinearStabilityOfExandingRicciSolitons}, it is shown that low dimensional solitons are stable and all solitons on 2-step nilpotent Lie groups are stable (among all metrics, not just homogeneous metrics).

\begin{openprob} Are non-compact, homogeneous Ricci solitons all (linearly) stable under the Ricci flow? 
\end{openprob}

If one restricts to a particular solvable Lie group which admits a soliton, then we know that the soliton is dynamically stable among homogeneous metrics from \cite{Bohm-Lafuente:TheRicciFlowOnSolvmanifoldsOfRealType}.  In general, i.e. among all metrics, stability is not well-understood.  See, also, \cite{Lott:DimReductionAndLongTimeBehaviorOfRicciFlow}.

\subsection{Maximal symmetry}\label{sec: maximal symmetry}
Given the special curvature properties of Einstein and Ricci soliton metrics on solvmanifolds, along with their uniqueness on a given Lie group, one might wonder what other special properties these spaces enjoy.

Consider a Lie group $S$ and a left-invariant metric $g$.  Recall, for $\phi\in Aut(S)$, $\phi^*g$ is another, isometric, left-invariant metric.  We say that $g$ has maximal symmetry if given an arbitrary left-invariant metric $g'$ on $S$ there is an automorphism $\phi \in Aut(G)$ such that
	$$Isom(S,g') \subset Isom(S,\phi^*g).$$

\begin{thm}\cite{GordonJablonski:EinsteinSolvmanifoldsHaveMaximalSymmetry} \label{thm: Einstein solv have maximal symmetry}
Einstein metrics on solvable Lie groups have maximal symmetry.
\end{thm}

The proof of this result relies on the pre-Einstein derivation, see Equation \ref{eqn: Ricci soliton defintion}.  A dynamical proof of the theorem above follows from the stability of Einstein solvmanifolds (among left-invariant metrics) under the Ricci flow.  While the above result does hold for unimodular solitons on solvmanifolds \cite{Jablo:MaximalSymmetryAndUnimodularSolvmanifolds}, it cannot hold for all non-unimodular solitons.  An infinitesimal version of it has been achieved \cite{GordonJablonski:RicciSolitonSolvmanifoldsHaveInfinitesimalMaximalSymmetry} where it is shown that one can get a containment of isometry algebras.

\begin{openprob} Is there a dynamical proof of infinitesimal symmetry for solvsolitons?
\end{openprob}

The technique of using Ricci flow that worked in the Einstein setting cannot work for solitons.  New tools are needed.

\begin{openprob} Are there algebraic criteria that determine when a given solvmanifold admits a maximal symmetry metric?
\end{openprob}

There is an interesting algebraic question that arises from studying maximal symmetry.  If $\mathfrak g$ is the isometry algebra of a completely solvable Lie algebra $\mathfrak s$, then decomposing $\mathfrak g$ into a Levi decomposition $\mathfrak g = \mathfrak g_1\ltimes \mathfrak g_2$ we are able to see that
	\begin{align*}
	\mathfrak s &= \mathfrak s_1\ltimes \mathfrak s_2
	\shortintertext{where}
	\mathfrak s_1 &\subset \mathfrak g_1 \mbox{ is the Iwasawa subgroup of } \mathfrak g_1
	\end{align*}
and the action of $\mathfrak s_1$ is the restriction of the adjoint action of $\mathfrak g_1$ on ideal $\mathfrak s_2$.  We call the decomposition above a pre-Levi decomposition (pre-Levi as it extends to a Levi decomposition of some isometry group).

The existence of infinitesimal maximal symmetry metrics then becomes an algebraic problem for completely solvable $\mathfrak s$ of finding `maximal pre-Levi decompositions'.

\begin{openprob}  Does every solvable Lie algebra admit a maximal pre-Levi decomposition?
\end{openprob}

For solvable with 2-step nilradicals, this has been considered in the forthcoming work \cite{EpsteinJablonski:MaximalSymmetryAndSolvmanifolds}.  There it is shown that these low step solvable Lie groups do admit metrics of infinitesimal maximal symmetry.

\subsection{Solitons as critical points of a geometric functional}
For compact homogeneous spaces, the Ricci flow is the negative gradient flow of the scalar curvature function and critical points of the scalar curvature function (restricted to the set of $G$-invariant, volume 1 metrics) are precisely the Einstein metrics.  This is not true for non-compact homogeneous spaces and another function has been proposed for study.

Motivated by GIT and the nilpotent setting (see Section \ref{subsec: GIT}), one may ask if homogeneous Ricci solitons on solvable groups are the maxima of the function
	$$F(g) = \frac{sc(g)^2}{tr~Ric_g^2}.$$
For unimodular solvmanifolds or those with codimension one nilradical, the answer is yes \cite{Lauret-Will:TheRicciPinchingFunctionalOnSolvmanifolds,Lauret-Will:TheRicciPinchingFunctionalOnSolvmanifolds2}!

\begin{openprob}  Are all Ricci solitons on solvable Lie groups maxima of the function $F(g) =\frac{sc(g)^2}{tr~Ric_g^2}$?
\end{openprob}

\bibliographystyle{amsalpha}

\begin{thebibliography}{HMW21}

\bibitem[ADF96]{Alekseevsky-Dotti-Ferraris:HomogeneosRicciPositive5manifolds}
D.~Alekseevsky, Isabel Dotti, and C.~Ferraris, \emph{Homogeneous {R}icci
  positive {$5$}-manifolds}, Pacific J. Math. \textbf{175} (1996), no.~1,
  1--12. \MR{1419469}

\bibitem[AK75]{AlekseevskiiKimelfeld:StructureOfHomogRiemSpacesWithZeroRicciCurv}
D.~V. Alekseevski{\u\i} and B.~N. Kimel{'}fel{'}d, \emph{Structure of
  homogeneous {R}iemannian spaces with zero {R}icci curvature}, Functional
  Anal. Appl. \textbf{9} (1975), no.~2, 97--102.

\bibitem[AL17]{Arroyo-Lafuente:TheAlekseevskiiConjectureInLowDimensions}
Romina~M. Arroyo and Ramiro~A. Lafuente, \emph{The {A}lekseevskii conjecture in
  low dimensions}, Math. Ann. \textbf{367} (2017), no.~1-2, 283--309.
  \MR{3606442}

\bibitem[And10]{Anderson:ASurveyOfEinsteinMetricsOn4Manifolds}
Michael~T. Anderson, \emph{A survey of {E}instein metrics on 4-manifolds},
  Handbook of geometric analysis, {N}o. 3, Adv. Lect. Math. (ALM), vol.~14,
  Int. Press, Somerville, MA, 2010, pp.~1--39. \MR{2743446 (2012a:53071)}

\bibitem[Arr11]{Arroyo:FiliformNilsolitonsOfDimension8}
Romina~M. Arroyo, \emph{Filiform nilsolitons of dimension 8}, Rocky Mountain J.
  Math. \textbf{41} (2011), no.~4, 1025--1043. \MR{2826522}

\bibitem[Arv15]{Arvanitoyeorgos:ProgressOnHomogeneousEinsteinManifoldsAndSomeOpenProblems}
Andreas Arvanitoyeorgos, \emph{Progress on homogeneous {E}instein manifolds and
  some open problems}, Bull. Greek Math. Soc. \textbf{58} (2010/15), 75--97.
  \MR{3585267}

\bibitem[ASS20a]{Arvanitoyeorgos-Sakane-Statha:HomogeneousEMetricsOnStiefelManifoldsAssociatedToFlagManifoldsWithTwoIsotropySummands}
Andreas Arvanitoyeorgos, Yusuke Sakane, and Marina Statha, \emph{Homogeneous
  {E}instein metrics on {S}tiefel manifolds associated to flag manifolds with
  two isotropy summands}, J. Symbolic Comput. \textbf{101} (2020), 189--201.
  \MR{4109715}

\bibitem[ASS20b]{Arvanitoyeorgos-Sakane-Statha:InvariantEinsteinMetricsOnSUnAndComplexStiefelManifolds}
\bysame, \emph{Invariant {E}instein metrics on {SU}(n)and complex {S}tiefel
  manifolds}, Tohoku Math. J. (2) \textbf{72} (2020), no.~2, 161--210.
  \MR{4116694}

\bibitem[B\"05]{Bohm:NonexistenceOfHomogeneousEinsteinMetrics}
Christoph B\"{o}hm, \emph{Non-existence of homogeneous {E}instein metrics},
  Comment. Math. Helv. \textbf{80} (2005), no.~1, 123--146. \MR{2130570}

\bibitem[BB78]{Berard-Bergery:SurLaCourbureDesMetriquesRiemanniennesInvariantesDesGroupesDeLieEtDesEspacesHomegenes}
Lionel B\'{e}rard-Bergery, \emph{Sur la courbure des m\'{e}triques
  riemanniennes invariantes des groupes de {L}ie et des espaces homog\`enes},
  Ann. Sci. \'{E}cole Norm. Sup. (4) \textbf{11} (1978), no.~4, 543--576.
  \MR{533067}

\bibitem[Ber21]{Berichon:TheAlekseevskiiConjectureIn9and10dimensions}
Rohin Berichon, \emph{The {A}lekseevskii conjecture in 9 and 10 dimensions},
  Differential Geom. Appl. \textbf{78} (2021), Paper No. 101782, 20.
  \MR{4280205}

\bibitem[Bes08]{Besse:EinsteinMflds}
Arthur~L. Besse, \emph{Einstein manifolds}, Classics in Mathematics,
  Springer-Verlag, Berlin, 2008, Reprint of the 1987 edition.

\bibitem[BK06]{Bohm-Kerr:LowDimensionalHomogeneousEinsteinManifolds}
Christoph B\"{o}hm and Megan~M. Kerr, \emph{Low-dimensional homogeneous
  {E}instein manifolds}, Trans. Amer. Math. Soc. \textbf{358} (2006), no.~4,
  1455--1468. \MR{2186982}

\bibitem[BL18]{Bohm-Lafuente:HomogeneousEinsteinMetricsOnEuclideanSpacesAreEinsteinSolvmanifolds}
Christoph B{\"o}hm and Ramiro Lafuente, \emph{Homogeneous {E}instein metrics on
  {E}uclidean spaces are {E}instein solvmanifolds}, arXiv:1811.12594 (2018).

\bibitem[BL19]{Bohm-Lafuente:TheRicciFlowOnSolvmanifoldsOfRealType}
Christoph B\"{o}hm and Ramiro~A. Lafuente, \emph{The {R}icci flow on
  solvmanifolds of real type}, Adv. Math. \textbf{352} (2019), 516--540.
  \MR{3964154}

\bibitem[BL21]{Bohm-Lafuente:NonCompactEinsteinManifoldsWithSymmetry}
Christoph B{\"o}hm and Ramiro Lafuente, \emph{Non-compact {E}instein manifolds
  with symmetry}, arXiv:2107.04210 (2021).

\bibitem[Boc46]{Bochner:VectorFieldsAndRicciCurvature}
S.~Bochner, \emph{Vector fields and {R}icci curvature}, Bull. Amer. Math. Soc.
  \textbf{52} (1946), 776--797. \MR{0018022 (8,230a)}

\bibitem[B{\"o}h04]{Bohm:HomogEinsteinMetricsAndSimplicialComplexes}
Christoph B{\"o}hm, \emph{Homogeneous {E}instein metrics and simplicial
  complexes}, J. Differential Geom. \textbf{67} (2004), no.~1, 79--165.

\bibitem[BWZ04]{Bohm-Wang-Ziller:AVariationalApproachforCompactHomogEinsteinMflds}
C.~B{\"o}hm, M.~Wang, and W.~Ziller, \emph{A variational approach for compact
  homogeneous {E}instein manifolds}, Geom. Funct. Anal. \textbf{14} (2004),
  no.~4, 681--733.

\bibitem[CCD18]{Chen-Chen-Deng:NonnaturallyReductiveEinsteinMetricsOnSOn}
Huibin Chen, Zhiqi Chen, and Shaoqiang Deng, \emph{Non-naturally reductive
  {E}instein metrics on {SO}(n)}, Manuscripta Math. \textbf{156} (2018),
  no.~1-2, 127--136. \MR{3783569}

\bibitem[CN19]{Chen-Nikonorov:InvariantEinsteinMetricsOnGeneralizedWallachSpaces}
Zhiqi Chen and Yuri\u{\i}~Gennadievich Nikonorov, \emph{Invariant {E}instein
  metrics on generalized {W}allach spaces}, Sci. China Math. \textbf{62}
  (2019), no.~3, 569--584. \MR{3905563}

\bibitem[CS21]{Chrysikos-Sakane:HomogeneousEinsteinMetricsOnNonKahlerCSpaces}
Ioannis Chrysikos and Yusuke Sakane, \emph{Homogeneous {E}instein metrics on
  non-{K}\"{a}hler {C}-spaces}, J. Geom. Phys. \textbf{160} (2021), Paper No.
  103996, 31. \MR{4176991}

\bibitem[DK08]{Dickinson-Kerr:TheGeometryOfCompactHomogeneousSpacesWithTwoIsotropySummands}
William Dickinson and Megan~M. Kerr, \emph{The geometry of compact homogeneous
  spaces with two isotropy summands}, Ann. Global Anal. Geom. \textbf{34}
  (2008), no.~4, 329--350. \MR{2447903}

\bibitem[DM82]{Dotti:RicciCurvUnimodularSolv}
Isabel Dotti~Miatello, \emph{Ricci curvature of left invariant metrics on
  solvable unimodular {L}ie groups}, Math. Z. \textbf{180} (1982), no.~2,
  257--263.

\bibitem[DZ79]{DatriZiller:NaturallyRedMetricsAndEinsteinMetricsOnCompactLieGroups}
J.~E. D'Atri and W.~Ziller, \emph{Naturally reductive metrics and {E}instein
  metrics on compact {L}ie groups}, Mem. Amer. Math. Soc. \textbf{18} (1979),
  no.~215, iii+72.

\bibitem[Ebe08]{Eberlein:prescribedRicciTensor}
Patrick Eberlein, \emph{Riemannian 2-step nilmanifolds with prescribed {R}icci
  tensor}, Geometric and probabilistic structures in dynamics, Contemp. Math.,
  vol. 469, Amer. Math. Soc., Providence, RI, 2008, pp.~167--195.

\bibitem[EJ21]{EpsteinJablonski:MaximalSymmetryAndSolvmanifolds}
Jonathan Epstein and Michael Jablonski, \emph{Maximal symmetry and
  solvmanifolds}, in progress (2021).

\bibitem[FC14]{Fernandez-Culma:ClassificationOfNilsolitonMetricsInDimensionSeven}
Edison~Alberto Fern\'{a}ndez-Culma, \emph{Classification of nilsoliton metrics
  in dimension seven}, J. Geom. Phys. \textbf{86} (2014), 164--179.
  \MR{3282320}

\bibitem[GJ19]{GordonJablonski:EinsteinSolvmanifoldsHaveMaximalSymmetry}
Carolyn Gordon and Michael Jablonski, \emph{Einstein solvmanifolds have maximal
  symmetry}, Journal of Differential Geometry \textbf{111} (2019), no.~1,
  1--38.

\bibitem[GJ21]{GordonJablonski:RicciSolitonSolvmanifoldsHaveInfinitesimalMaximalSymmetry}
\bysame, \emph{Ricci soliton solvmanifolds have infinitesimal maximal
  symmetry}, preprint (2021).

\bibitem[Gor80]{Gordon:RiemannianIsometryGroupsContainingTransitiveReductiveSubgroups}
Carolyn Gordon, \emph{Riemannian isometry groups containing transitive
  reductive subgroups}, Math. Ann. \textbf{248} (1980), no.~2, 185--192.
  \MR{573347 (81e:53030)}

\bibitem[GP10]{GlickensteinPayne:RicciFlow3dimUnimodLieAlgs}
David Glickenstein and Tracy~L. Payne, \emph{Ricci flow on three-dimensional,
  unimodular metric {L}ie algebras}, Comm. Anal. Geom. \textbf{18} (2010),
  no.~5, 927--961. \MR{2805148 (2012k:53129)}

\bibitem[GW85]{GordonWilson:TheFineStructureOfTransitiveRiemannianIsometryGroups}
Carolyn~S. Gordon and Edward~N. Wilson, \emph{The fine structure of transitive
  {R}iemannian isometry groups. {I}}, Trans. Amer. Math. Soc. \textbf{289}
  (1985), no.~1, 367--380. \MR{779070 (86g:53056)}

\bibitem[GW88]{GordonWilson:IsomGrpsOfRiemSolv}
\bysame, \emph{Isometry groups of {R}iemannian solvmanifolds}, Trans. Amer.
  Math. Soc. \textbf{307} (1988), no.~1, 245--269.

\bibitem[He12]{He:CohomogeneityOneManifoldsWithASmallFamilyOfInvariantMetrics}
Chenxu He, \emph{Cohomogeneity one manifolds with a small family of invariant
  metrics}, Geom. Dedicata \textbf{157} (2012), 41--90. \MR{2893479}

\bibitem[Heb98]{Heber}
Jens Heber, \emph{Noncompact homogeneous {E}instein spaces}, Invent. Math.
  \textbf{133} (1998), no.~2, 279--352.

\bibitem[HMW21]{Hall-Murphy-Waldron:CompactHermitianSymmetricSpaces}
Stuart~James Hall, Thomas Murphy, and James Waldron, \emph{Compact {H}ermitian
  symmetric spaces, coadjoint orbits, and the dynamical stability of the
  {R}icci flow}, J. Geom. Anal. \textbf{31} (2021), no.~6, 6195--6218.
  \MR{4267642}

\bibitem[Jab08]{Jablo:Thesis}
Michael Jablonski, \emph{Real geometric invariant theory and {R}icci soliton
  metrics on two-step nilmanifolds}, Thesis (May 2008).

\bibitem[Jab10]{Jablo:DetectingOrbitsAlongSubvarietiesViaTheMomentMap}
\bysame, \emph{Detecting orbits along subvarieties via the moment map},
  M\"unster J. Math. \textbf{3} (2010), 67--88.

\bibitem[Jab11a]{Jablo:ConceringExistenceOfEinstein}
\bysame, \emph{Concerning the existence of {E}instein and {R}icci soliton
  metrics on solvable lie groups}, Geometry \& Topology \textbf{15} (2011),
  no.~2, 735--764.

\bibitem[Jab11b]{Jablo:ModuliOfEinsteinAndNoneinstein}
\bysame, \emph{Moduli of {E}instein and non-{E}instein nilradicals}, Geom.
  Dedicata \textbf{152} (2011), no.~1, 63--84.

\bibitem[Jab14]{Jablo:HomogeneousRicciSolitonsAreAlgebraic}
\bysame, \emph{Homogeneous {R}icci solitons are algebraic}, Geom. Topol.
  \textbf{18} (2014), no.~4, 2477--2486. \MR{3268781}

\bibitem[Jab15a]{Jablo:HomogeneousRicciSolitons}
\bysame, \emph{Homogeneous {R}icci solitons}, J. Reine Angew. Math.
  \textbf{699} (2015), 159--182. \MR{3305924}

\bibitem[Jab15b]{Jablo:StronglySolvable}
\bysame, \emph{Strongly solvable spaces}, Duke Math. J. \textbf{164} (2015),
  no.~2, 361--402. \MR{3306558}

\bibitem[Jab19]{Jablo:MaximalSymmetryAndUnimodularSolvmanifolds}
\bysame, \emph{Maximal symmetry and unimodular solvmanifolds}, Pacific J. Math.
  \textbf{298} (2019), no.~2, 417--427. \MR{3936023}

\bibitem[Jen69]{Jensen:HomogEinsteinSpacesofDim4}
Gary~R. Jensen, \emph{Homogeneous {E}instein spaces of dimension four}, J.
  Differential Geometry \textbf{3} (1969), 309--349.

\bibitem[Jen71]{Jensen:TheScalarCurvatureOfLeftInvariantRiemannianMetrics}
\bysame, \emph{The scalar curvature of left-invariant {R}iemannian metrics},
  Indiana Univ. Math. J. \textbf{20} (1970/1971), 1125--1144. \MR{0289726 (44
  \#6914)}

\bibitem[JP17]{JP:TowardsTheAlekseevskiiConjecture}
Michael Jablonski and Peter Petersen, \emph{A step towards the {A}lekseevskii
  conjecture}, Math. Ann. \textbf{368} (2017), no.~1-2, 197--212. \MR{3651571}

\bibitem[JPW14]{JPW2:OnTheLinearStabilityOfExandingRicciSolitons}
Michael Jablonski, Peter Petersen, and Michael~Bradford Williams, \emph{On the
  linear stability of expanding {R}icci solitons}, arXiv:1409.3251 [math.DG]
  (2014).

\bibitem[JPW16]{JPW:LinearStabilityOfAlgebraicRicciSolitons}
\bysame, \emph{Linear stability of algebraic {R}icci solitons}, J. Reine Angew.
  Math. \textbf{713} (2016), 181--224. \MR{3483629}

\bibitem[Kir84]{Kirwan}
Frances~Clare Kirwan, \emph{Cohomology of quotients in symplectic and algebraic
  geometry}, Mathematical Notes 31, Princeton University Press, Princeton, New
  Jersey, 1984.

\bibitem[Kot10]{Kotschwar:BackwardsUniquenessRF}
B.L. Kotschwar, \emph{Backwards uniqueness for the {R}icci flow}, Int. Math.
  Res. Not. IMRN (2010), no.~21, 4064--4097.

\bibitem[Kr{\"o}20]{Kroncke:StabilityOfEinsteinMetricsUnderRicciFlow}
Klaus Kr{\"o}ncke, \emph{Stability of {E}instein metrics under {R}icci flow},
  Comm. Anal. Geom. \textbf{28} (2020), no.~2, 351--394. \MR{4101342}

\bibitem[Lau01a]{Lauret:RicciSolitonHomogeneousNilmanifolds}
Jorge Lauret, \emph{Ricci soliton homogeneous nilmanifolds}, Math. Ann.
  \textbf{319} (2001), no.~4, 715--733.

\bibitem[Lau01b]{Lauret:StandardEinsteinSolvAsCriticalPoints}
\bysame, \emph{Standard {E}instein solvmanifolds as critical points}, Q. J.
  Math. \textbf{52} (2001), no.~4, 463--470.

\bibitem[Lau02]{Lauret:FindingEinsteinSolvmanifoldsByAVariantialMethod}
\bysame, \emph{Finding {E}instein solvmanifolds by a variational method}, Math.
  Z. \textbf{241} (2002), no.~1, 83--99. \MR{1930986}

\bibitem[Lau03a]{Lauret:DegenerationsOfLieAlgebrasAndGeometryOfLieGroups}
\bysame, \emph{Degenerations of {L}ie algebras and geometry of {L}ie groups},
  Differential Geom. Appl. \textbf{18} (2003), no.~2, 177--194. \MR{1958155}

\bibitem[Lau03b]{Lauret:MomentMapVarietyLieAlgebras}
\bysame, \emph{On the moment map for the variety of {L}ie algebras}, J. Funct.
  Anal. \textbf{202} (2003), no.~2, 392--423.

\bibitem[Lau09]{Lauret:EinsteinSolvandNilsolitonsCordobaConf2007}
\bysame, \emph{Einstein solvmanifolds and nilsolitons}, New developments in
  {L}ie theory and geometry, Contemp. Math., vol. 491, Amer. Math. Soc.,
  Providence, RI, 2009, pp.~1--35.

\bibitem[Lau10]{Lauret:EinsteinSolvmanifoldsAreStandard}
\bysame, \emph{Einstein solvmanifolds are standard}, Ann. of Math. (2)
  \textbf{172} (2010), no.~3, 1859--1877.

\bibitem[Lau11a]{Lauret:RicciFlowForSimplyConnectedNilmanifolds}
\bysame, \emph{The {R}icci flow for simply connected nilmanifolds}, Comm. Anal.
  Geom. \textbf{19} (2011), no.~5, 831--854. \MR{2886709}

\bibitem[Lau11b]{Lauret:SolSolitons}
\bysame, \emph{Ricci soliton solvmanifolds}, J. Reine Angew. Math. \textbf{650}
  (2011), 1--21.

\bibitem[Lau13]{Lauret:RicciFlowOfHomogeneousManifoldsAndItsSolitons}
\bysame, \emph{Ricci flow of homogeneous manifolds}, Math. Z. \textbf{274}
  (2013), no.~1-2, 373--403. \MR{3054335}

\bibitem[Lau21]{Lauret:OnTheStabilityOfHomogeneousEinsteinManifolds}
\bysame, \emph{On the stability of homogeneous {E}instein manifolds},
  arXiv:2105.06336 (2021).

\bibitem[LL14]{LauretLafuente:StructureOfHomogeneousRicciSolitonsAndTheAlekseevskiiConjecture}
Ramiro Lafuente and Jorge Lauret, \emph{Structure of homogeneous {R}icci
  solitons and the {A}lekseevskii conjecture}, J. Differential Geom.
  \textbf{98} (2014), no.~2, 315--347. \MR{3263520}

\bibitem[LO14]{Lauret-Oscari:OnNonsingular2stepNilpotentLieAlgebras}
Jorge Lauret and David Oscari, \emph{On non-singular 2-step nilpotent {L}ie
  algebras}, Math. Res. Lett. \textbf{21} (2014), no.~3, 553--583. \MR{3272030}

\bibitem[Loh94]{Lohkamp:MetricsOfNegativeRicciCurvature}
Joachim Lohkamp, \emph{Metrics of negative {R}icci curvature}, Ann. of Math.
  (2) \textbf{140} (1994), no.~3, 655--683. \MR{1307899}

\bibitem[Lot10]{Lott:DimReductionAndLongTimeBehaviorOfRicciFlow}
John Lott, \emph{Dimensional reduction and the long-time behavior of {R}icci
  flow}, Comment. Math. Helv. \textbf{85} (2010), no.~3, 485--534.

\bibitem[LW11]{LauretWill:EinsteinSolvExistandNonexist}
Jorge Lauret and Cynthia Will, \emph{Einstein solvmanifolds: existence and
  non-existence questions}, Math. Ann. \textbf{350} (2011), no.~1, 199--225.
  \MR{2785768 (2012i:53040)}

\bibitem[LW19]{Lauret-Will:TheRicciPinchingFunctionalOnSolvmanifolds}
Jorge Lauret and Cynthia~E. Will, \emph{The {R}icci pinching functional on
  solvmanifolds}, Q. J. Math. \textbf{70} (2019), no.~4, 1281--1304.
  \MR{4045101}

\bibitem[LW20]{Lauret-Will:TheRicciPinchingFunctionalOnSolvmanifolds2}
\bysame, \emph{The {R}icci pinching functional on solvmanifolds {II}}, Proc.
  Amer. Math. Soc. \textbf{148} (2020), no.~6, 2601--2607. \MR{4080900}

\bibitem[Man66]{Manturov:HomogeneousRiemannianSpacesWwithAnIrreducibleRotationGroup}
O.~V. Manturov, \emph{Homogeneous {R}iemannian spaces with an irreducible
  rotation group}, Trudy Sem. Vektor. Tenzor. Anal. \textbf{13} (1966),
  68--145. \MR{0210031}

\bibitem[Nab07]{Naber:NoncompactShrinking4SolitonsWithNonnegativeCurvature}
Aaron Naber, \emph{Noncompact shrinking 4-solitons with nonnegative curvature},
  arXiv:0710.5579 [math.DG] (2007).

\bibitem[Nik98]{Nikonorov:TheScalarCurvFuctionalAndHomogEinsteinMetricsOnLieGroups}
Yu.~G. Nikonorov, \emph{The scalar curvature functional and homogeneous
  {E}instein metrics on {L}ie groups}, Sibirsk. Mat. Zh. \textbf{39} (1998),
  no.~3, 504--509, ii. \MR{1639512 (99h:53061)}

\bibitem[Nik04]{Nikonorv:CompactHomogeneousEinstein7manifolds}
\bysame, \emph{Compact homogeneous {E}instein 7-manifolds}, Geom. Dedicata
  \textbf{109} (2004), 7--30. \MR{2113184}

\bibitem[Nik05]{Nikonorov:NoncompactHomogEinstein5manifolds}
\bysame, \emph{Noncompact homogeneous {E}instein 5-manifolds}, Geom. Dedicata
  \textbf{113} (2005), 107--143.

\bibitem[Nik06]{Nikitenko:Seven-dimensionalHomogeneousEinsteinManifoldsOfNegativeSectionalCurvature}
E.~V. Nikitenko, \emph{Seven-dimensional homogeneous {E}instein manifolds of
  negative sectional curvature}, Mat. Tr. \textbf{9} (2006), no.~1, 101--116.
  \MR{2251332 (2007d:53078)}

\bibitem[Nik07]{Nikonorov:OnEinsteinExtensionsOfNilpotentMetricLieAlgebras}
Yu.~G. Nikonorov, \emph{On {E}instein extensions of nilpotent metric {L}ie
  algebras}, Mat. Tr. \textbf{10} (2007), no.~1, 164--190. \MR{2485371}

\bibitem[Nik11]{Nikolayevsky:EinsteinSolvmanifoldsandPreEinsteinDerivation}
Y.~Nikolayevsky, \emph{Einstein solvmanifolds and the pre-{E}instein
  derivation}, Trans. Amer. Math. Soc. \textbf{363} (2011), 3935--3958.

\bibitem[Nik16]{Nikonorov:ClassificationOfGeneralizedWallachSpaces}
Yu.~G. Nikonorov, \emph{Classification of generalized {W}allach spaces}, Geom.
  Dedicata \textbf{181} (2016), 193--212. \MR{3475745}

\bibitem[NN05]{NikitenkoNikonorov:Six-dimensionalEinsteinSolvmanifolds}
E.~V. Nikitenko and Yu.~G. Nikonorov, \emph{Six-dimensional {E}instein
  solvmanifolds}, Mat. Tr. \textbf{8} (2005), no.~1, 71--121. \MR{1955023
  (2006a:53043)}

\bibitem[NR03]{Nikonorov-Rodionov:CompactHomogeneousEinstein6manifolds}
Yu.~G. Nikonorov and E.~D. Rodionov, \emph{Compact homogeneous {E}instein
  6-manifolds}, Differential Geom. Appl. \textbf{19} (2003), no.~3, 369--378.
  \MR{2013101}

\bibitem[Osc14]{Oscari:OnTheExistenceOfNilsolitonsOn2stepNilpotentLieGroups}
David Oscari, \emph{On the existence of nilsolitons on 2-step nilpotent {L}ie
  groups}, Adv. Geom. \textbf{14} (2014), no.~3, 483--497. \MR{3228895}

\bibitem[Pay10]{Payne:ExistenceofSolitononNil}
Tracy~L. Payne, \emph{The existence of soliton metrics for nilpotent {L}ie
  groups}, Geom. Dedicata \textbf{145} (2010), 71--88.

\bibitem[Pay12]{Payne:GeometricInvariantsForNilpotentMetricLieAlgebrasWithApplicationsToModuliSpacesOfNilsolitonMetrics}
\bysame, \emph{Geometric invariants for nilpotent metric {L}ie algebras with
  applications to moduli spaces of nilsoliton metrics}, Ann. Global Anal. Geom.
  \textbf{41} (2012), no.~2, 139--160. \MR{2876692}

\bibitem[PW09]{Petersen-Wylie:OnGradientRicciSolitonsWithSymmetry}
Peter Petersen and William Wylie, \emph{On gradient {R}icci solitons with
  symmetry}, Proc. Amer. Math. Soc. \textbf{137} (2009), no.~6, 2085--2092.
  \MR{2480290 (2010a:53073)}

\bibitem[Wan12]{Wang:EinsteinMetricsFromSymmetryAndBundleConstractionsASurvey}
McKenzie Y.-K. Wang, \emph{Einstein metrics from symmetry and bundle
  constructions: a sequel}, Differential geometry, Adv. Lect. Math. (ALM),
  vol.~22, Int. Press, Somerville, MA, 2012, pp.~253--309. \MR{3076055}

\bibitem[Wil03]{Will:Rank1EinsteinSolvOfDim7}
C.E. Will, \emph{Rank-one einstein solvmanifolds of dimension 7}, Diff. Geom.
  Appl. \textbf{19} (2003), 307--318.

\bibitem[Wil10]{Will:CurveOfNonEinsteinNilradicals}
Cynthia Will, \emph{A curve of nilpotent {L}ie algebras which are not
  {E}instein nilradicals}, Monatsh. Math. \textbf{159} (2010), no.~4, 425--437.

\bibitem[Wil11]{Will:TheSpaceOfSolsolitonsInLowDimensions}
\bysame, \emph{The space of solvsolitons in low dimensions}, Ann. Global Anal.
  Geom. \textbf{40} (2011), no.~3, 291--309. \MR{2831460 (2012j:53053)}

\bibitem[Wol68]{Wolf:TheGeometryAndStructureOfIsotropyIrreducibleHomogeneousSpaces}
Joseph~A. Wolf, \emph{The goemetry and structure of isotropy irreducible
  homogeneous spaces}, Acta Math. \textbf{120} (1968), 59--148. \MR{223501}

\bibitem[WW16]{WilliamsWu:DynamicalStabilityOfAlgebraicRicciSolitons}
Michael~Bradford Williams and Haotian Wu, \emph{Dynamical stability of
  algebraic {R}icci solitons}, J. Reine Angew. Math. \textbf{713} (2016),
  225--243. \MR{3483630}

\bibitem[WZ86]{Wang-Ziller:ExistenceAndNonexistenceOfHomogEinstein}
McKenzie~Y. Wang and Wolfgang Ziller, \emph{Existence and nonexistence of
  homogeneous {E}instein metrics}, Invent. Math. \textbf{84} (1986), no.~1,
  177--194.

\bibitem[WZ90]{Wang-Ziller:EinsteinMetricsOnPrincipalTorusBundles}
\bysame, \emph{Einstein metrics on principal torus bundles}, J. Differential
  Geom. \textbf{31} (1990), no.~1, 215--248. \MR{1030671}

\bibitem[YD17]{Yan-Deng:EinsteinMetricsOnCompactSimpleLieGroupsAttachedToStandardTriples}
Zaili Yan and Shaoqiang Deng, \emph{Einstein metrics on compact simple {L}ie
  groups attached to standard triples}, Trans. Amer. Math. Soc. \textbf{369}
  (2017), no.~12, 8587--8605. \MR{3710636}

\bibitem[YD18]{Yan-Deng:InvariantEinsteinMetricsOnCertainCompactSemisimple}
\bysame, \emph{Invariant {E}instein metrics on certain compact semisimple {L}ie
  groups}, Differential Geom. Appl. \textbf{59} (2018), 138--153. \MR{3804826}

\end{thebibliography}

\providecommand{\bysame}{\leavevmode\hbox to3em{\hrulefill}\thinspace}
\providecommand{\MR}{\relax\ifhmode\unskip\space\fi MR }
\providecommand{\MRhref}[2]{%
  \href{http://www.ams.org/mathscinet-getitem?mr=#1}{#2}
}
\providecommand{\href}[2]{#2}

\end{document}